\documentclass[fleqn]{article}

\newcommand{\R}{\mathbb{R}}

\newcommand{\bq}{\mbox{{\boldmath $\mathcal{Q}$}}}
\newcommand{\bp}{\mbox{{\boldmath $\mathcal{P}$}}}
\newcommand{\eps}{\varepsilon}

\newtheorem{thm}{Theorem}[section]
\newtheorem{lem}[thm]{Lemma}
\newtheorem{prop}[thm]{Proposition}
\newtheorem{cor}[thm]{Corollary}

\newtheorem{remark}[thm]{Remark}
\usepackage{amsmath}
\usepackage{amssymb}

\begin{document}
\title{Nonnegative Matrix Factorization and I-Divergence Alternating
Minimization}
\author{Lorenzo Finesso\footnote{The authors have been supported in part by
the European Community's Human Potential Programme under contract
HPRN-CT-2000-00100, DYNSTOCH.} \\ ISIB--CNR \\ Corso Stati Uniti, 4 \\ 35127
Padova -- Italy \\ {\tt finesso@isib.cnr.it}
\and Peter Spreij\footnote{corresponding author} \\ Korteweg-de Vries
Institute for
Mathematics \\
Universiteit van Amsterdam \\
Plantage Muidergracht 24 \\
1018 TV Amsterdam -- The Netherlands \\
Email: {\tt spreij@science.uva.nl} \\
Phone: +31-20-5256070 \\
Fax: +31-20-5255101 \date{}} \maketitle
\begin{abstract}
\noindent In this paper we consider the Nonnegative Matrix
Factorization (NMF) problem: given an (elementwise) nonnegative
matrix $V \in \R_+^{m\times n}$ find, for assigned $k$,
nonnegative matrices $W\in\R_+^{m\times k}$ and $H\in\R_+^{k\times
n}$ such that $V=WH$. Exact, non trivial, nonnegative
factorizations do not always exist, hence it is interesting to
pose the approximate NMF problem. The criterion which is commonly
employed is I-divergence between nonnegative matrices. The problem
becomes that of finding, for assigned $k$, the factorization $WH$
closest to $V$ in I-divergence. An i\-te\-ra\-ti\-ve algorithm, EM
like, for the construction of the best pair $(W, H)$ has been
proposed in the literature. In this paper we interpret the
algorithm as an alternating minimization procedure \`a la
Csisz\'ar-Tusn\'ady and investigate some of its stability
properties. NMF is widespreading as a data analysis method in
applications for which the positivity constraint is relevant.
There are other data analysis methods which impose some form of
nonnegativity: we discuss here the connections between NMF and
Archetypal Analysis.

\end{abstract}

\newpage

\section{Introduction}

The approximate Nonnegative Matrix Factorization (NMF) of
nonnegative matrices is a data analysis technique only recently
introduced \cite{leeseung1999, sullivan2000}. Roughly speaking the
problem is to find, for a given nonnegative matrix $V \in
\R_+^{m\times n}$, and an assigned $k$, a pair of nonnegative
matrices $W\in\R_+^{m\times k}$ and $H\in\R_+^{k\times n}$ such
that, in an appropriate sense, $V \approx WH$. In
\cite{leeseung1999} EM like algorithms for the construction of a
factorization have been proposed. The algorithms have been later
derived in \cite{leeseung2001} by using an {\em ad-hoc} auxiliary
function, a common approach in deriving EM algorithms. In
\cite{sullivan2000} the connection with the classic alternating
minimization of the I-divergence \cite{ct1984} has been pointed
out but not fully investigated. In this paper we pose the NMF
problem as a minimum I-divergence problem that can be solved by
alternating minimization and derive, from this point of view, the
algorithm proposed in \cite{leeseung1999}. There are alternative
approaches to approximate nonnegative matrix factorization. For
instance, recently, see~\cite{laa2004}, results have been obtained
for the approximate factorization (w.r.t. the Frobenius norm) of
{\em symmetric} nonnegative matrices.

Although only recently introduced the NMF has found many
applications as a data reduction procedure and has been advocated
as an alternative to Principal Components Analysis (PCA) in cases
where the positivity constraint is relevant (typically image
analysis). The title of \cite{sullivan2000} is a clear indication
of this point of view, but a complete analysis of the relations
between NMF and PCA is still lacking.
Our interest in NMF stems from the system theoretic problem of
approximate realization (or order reduction) of Hidden Markov
Models. Partial results have already been obtained \cite{fs2002}.

This paper is organized as follows. In section~\ref{sec:problem}
we pose the approximate nonnegative matrix factorization problem,
define the I-divergence between matrices and discuss the solution
proposed in \cite{leeseung1999, leeseung2001}. In
section~\ref{section:lift} we pave the way for the alternating
minimization algorithm presenting the properly lifted version of
the minimization problem and solving the two partial minimizations
in the style of Csisz\'ar and Tusn\'ady \cite{ct1984}. In
section~\ref{section:altmin} we construct the alternating
minimization algorithm and compute the iteration gain. One of the
advantages of working with the lifted problem is that it sheds a
new light also on the derivation of the algorithm via auxiliary
functions given in \cite{leeseung2001}. In
section~\ref{section:auxfunc} we will use the results of
section~\ref{section:lift} to construct a very natural auxiliary
function to solve the original problem. A discussion of the
convergence properties of the algorithm is given in
section~\ref{section:convprop}. In the concluding
section~\ref{section:otherprob} we establish a connection between
the approximate NMF problem and the Archetypal Analysis algorithm
of Cutler and Breiman \cite{cb1994}. The present paper is an
extended version of~\cite{MTNS2004}.

\section{Preliminaries and problem statement}\label{sec:problem}

The NMF is a long standing problem in linear algebra \cite{haz,
phs}. It can be stated as follows. Given $V \in \R_+^{m\times n}$,
and $1 \le k \le \min\{m,n\}$, find a pair of matrices
$W\in\R_+^{m\times k}$ and $H\in\R_+^{k\times n}$ such that $V =
WH$.  The smallest $k$ for which a factorization exists is called
the positive rank of $V$, denoted ${\rm prank}(V)$. This
definition implies that ${\rm rank}(V) \le {\rm prank}(V) \le
\min\{m,n\}$. It is well known that ${\rm prank}(V)$ can assume
all intermediate values, depending on $V$. Examples for which
nonnegative factorizations do not exist, and examples for which
factorization is possible only for $k > {\rm rank}(V)$ have been
constructed in the literature \cite{haz}. The ${\rm prank}$ has
been characterized only for special classes of matrices \cite{phs}
and algorithms for the construction of a NMF of a general positive
matrix are not known.

The \emph{approximate} NMF has been recently introduced in
\cite{leeseung1999} independently from the exact NMF problem. The
set-up is the same, but instead of exact factorization it is
required that $V \approx WH$ in an appropriate sense. In
\cite{leeseung1999}, and in this paper, the approximation is to be
understood in the sense of minimum I-divergence. For two
nonnegative numbers $p$ and $q$ the I-divergence is defined as
\begin{equation*}
D(p||q)= p\log\frac{p}{q}-p+q,
\end{equation*}
with the conventions $0/0=0$, $0\log 0=0$ and $p/0=\infty$ for
$p>0$. From the inequality $x\log x\geq x-1$ it follows that
$D(p||q)\geq 0$ with equality iff $p=q$. For two nonnegative
matrices $M=(M_{ij})$ and $N=(N_{ij})$, of the same size, the
I-divergence is defined as
\begin{equation*}
D(M||N)= \sum_{ij}D(M_{ij}||N_{ij}).
\end{equation*}
Again it follows that $D(M||N)\geq 0$ with equality iff $M=N$. For
nonnegative vectors or tensors of the same size a similar
definition applies.

\noindent The problem of approximate NMF is to find for given $V$
and a fixed number $k$ (often referred to as the {\em inner size}
of the factorization)
\begin{equation}\label{eq:dmin}
\arg \min_{W,H}D(V||WH).
\end{equation}
The function $D: (W,H) \to D(V||WH)$ will sometimes be referred to
as the {\em objective function}. The {\em domain} of $D$ is the
set of pairs $(W,H)$ with nonnegative entries. The {\em interior}
of the domain is the subset of pairs $(W,H)$ with positive ($>0$)
entries, whereas pairs on the {\em boundary} have at least one
entry equal to zero.

Although the objective function $(W,H)\mapsto D(V||WH)$ is easily
seen to be convex in $W$ and $H$ separately, it is not jointly
convex in the two variables. Hence $(W,H)\mapsto D(V||WH)$ may
have several (local) minima and saddle points, that may prevent
numerical minimization algorithms to converge to the global
minimizer. However $D(V||WH)$ cannot have a local maximum in an
interior point $(W_0,H_0)$, because then also $W\mapsto
D(V||WH_0)$ would have a local maximum in $W_0$, which contradicts
convexity. Local maxima at the boundary are not a priori excluded.

\medskip
It is not immediately obvious that the approximate NMF problem
admits a solution. The following result is therefore relevant.
\begin{prop}\label{prop:exist}
The minimization problem~(\ref{eq:dmin}) has a solution.
\end{prop}
The proof of this proposition is deferred to
section~\ref{section:altmin}.

\medskip \noindent Notice that, increasing the inner size from $k$ to
$k+1$, the optimal value of the objective function decreases. This
follows from the fact that one can trivially embed the
factorization problem with inner size $k$ into the problem with
inner size $k+1$ simply adding a zero last column to the optimal
$W$ and an arbitrary last row to the optimal $H$ of the problem
with inner size $k$. Unfortunately, unlike the SVD of a matrix,
the best approximations with increasing $k$ are not embedded one
into another. For increasing $k$ the computations are to be
carried out anew.

Although, according to proposition~\ref{prop:exist}, a solution to
the minimization problem exists, it will certainly not be unique.
In order to rule out too many trivial multiple solutions, we
impose the condition that $H$ is row stochastic, so
$\sum_jH_{lj}=1$ for all $l$. This is not a restriction. Indeed,
first we exclude without loss of generality the case where $H$ has
one or more zero rows, since we would then in fact try to minimize
the I-divergence with inner size smaller than $k$. Let $h$ be the
diagonal matrix with elements $h_i=\sum_j H_{ij}$, then
$WH=\tilde{W}\tilde{H}$ with $\tilde{W}=Wh$, $\tilde{H}=h^{-1}H$
and $\tilde{H}$ is by construction row stochastic. The convention
that $H$ is row stochastic still does not rule out non-uniqueness.
Think e.g.\ of post-multiplying $W$ with a permutation matrix
$\Pi$ and pre-multiplying $H$ with $\Pi^{-1}$.

Let $e_n$ ($e_n^\top$) be the column (row) vector of size $n$
whose elements are all equal to one. Given $k$, the (constrained)
problem we will look at from now on is
\begin{equation}\label{maxc}
\min_{W,H: He_m=e_k} D(V||WH).
\end{equation}
For the sake of brevity we will often write $e$ for a vector of
$1$'s of generic size. The constraint in the previous problem will
then read as $He=e$.

\medskip
To carry out the minimization numerically, Lee and
Seung~\cite{leeseung1999,leeseung2001} proposed the following
iterative algorithm. Denoting by $W^t$ and $H^t$ the matrices at
step $t$, the update equations are
\begin{align}
W^{t+1}_{il} & =  W^t_{il} \sum_j
\frac{H^t_{lj}V_{ij}}{(W^tH^t)_{ij}}\label{eq:wbar2}\\
H^{t+1}_{lj} & =  H^t_{lj} \sum_i
\frac{W^t_{il}V_{ij}}{(W^tH^t)_{ij}} \Big/ \sum_{ij}
\frac{W^t_{il}H^t_{lj}V_{ij}}{(W^tH^t)_{ij}}.\label{eq:hbar2}
\end{align}
The initial condition $(W^0, H^0)$ will always be assumed to be in
the interior of the domain. Only a partial justification for this
algorithm is given in \cite{leeseung2001}, although the update
steps~(\ref{eq:wbar2}) and~(\ref{eq:hbar2}) are like those in the
EM algorithm, known from statistics, see~\cite{em}. Likewise the
convergence properties of the algorithm are unclear. In the next
section the minimization problem will be cast in a different way
to provide more insight in the specific form of the update
equations and on the convergence properties of the algorithm.

\medskip
We will now show that the $V$ matrix in the approximate NMF
problem can always be taken as a probability matrix $P$ i.e. such
that $P_{ij}\ge 0, \sum_{ij}P_{ij}=1$. This will pave the way for
the probabilistic interpretation of the exact and approximate NMF
problems to be given later.

Let $P =\frac{1}{e^\top Ve}V$, $Q_- =\frac{1}{e^\top We}W$,
$w=e^\top We$ and $Q_+ = H$. Notice that $e^\top Pe=e^\top Q_-e=1$
and $Q_+e=e$. Using the definition of divergence and elementary
computations, we obtain the decomposition
\[
D(V||WH)=e^\top Ve\, D(P||Q_-Q_+)+ D(e^\top Ve||w).
\]
Hence, since the number $e^\top V e$ is known, minimizing
$D(V||WH)$ w.r.t. $(W,H)$ is equivalent to minimizing
$D(P||Q_-Q_+)$ w.r.t.\ $(Q_-, Q_+)$ and $D(e^\top Ve||w)$ w.r.t.\
$w$. The minimizers of the three problems satisfy the relations
$W^*=e^\top Ve\, Q_-^*$, $H^* = Q_+^*$, and $w^*=e^\top Ve$.
Minimizing $D(V||WH)$ is therefore equivalent to minimizing
$D(P||Q_-Q_+)$. This enables us to give the problem a
probabilistic interpretation. Indeed,
\begin{equation}\label{eq:divnew}
D(P||Q_-Q_+)=\sum_{ij} D(P_{ij}||(Q_-Q_+)_{ij}) =
\sum_{ij}P_{ij}\log\frac{P_{ij}}{(Q_-Q_+)_{ij}},
\end{equation}
which is the usual I-divergence (Kullback-Leibler distance)
between (finite) probability measures. This will be exploited in
later sections.  From now on we will always consider the following
problem. Given the probability matrix $P$ and the integer $k$ find
\[
\min_{Q_-, Q_+ : Q_+e=e} D(P||Q_-Q_+).
\]
For typographical reasons we often, but not always, denote the
entries of $P$ by $P(ij)$ instead of $P_{ij}$ and likewise for
other matrices. \medskip\\
The minimization algorithm is easily seen to be {\em invariant
under the previous normalizations}. Let $Q_-^t=\frac{W^t}{e^\top
W^te}$ and $Q_-^t=H^t$. Substitute the definitions of $(P, Q_-^t,
Q_+^t)$ into~(\ref{eq:wbar2}) and~(\ref{eq:hbar2}) and use the
easily verified fact that $e^\top W^te=e^\top Ve$ for $t\geq 1$ to
obtain the update equations in the new notations
\begin{align}
Q_-^{t+1}(il) & = Q_-^t(il) \sum_j
\frac{Q_+^t(lj)P(ij)}{(Q_-^tQ_+^t)(ij)}\label{eq:q-bar2}\\
Q_+^{t+1}(lj) & = Q_+^t(lj) \sum_i
\frac{Q_-^t(il)P(ij)}{(Q_-^tQ_+^t)(ij)} \Big/ \sum_{ij}
\frac{Q_-^t(il)Q_+^t(lj)P(ij)}{(Q_-^tQ_+^t)(ij)}.\label{eq:q+bar2}
\end{align}
\rm

\section{Lifted version of the problem}\label{section:lift}

In this section we lift the I-divergence minimization problem to
an equivalent minimization problem where the `matrices' (we should
speak of {\em tensors}) have three indices.

\subsection{Setup}
Let be given a probability matrix $P$ (i.e. $P(ij) \ge 0, \,\,
\sum_{ij}P(ij)=1$) and an integer $k\le \min\{m,n\}$. We introduce
the following sets
\begin{align*}
\bp  = & \left\{\mathbf{P}\in\mathbb{R}^{m\times k\times n}_+ \,:
\,\, \sum_l\mathbf{P}(ilj)=P(ij)\right\}, \\ \\
\bq  = & \big\{\mathbf{Q}\in\mathbb{R}^{m\times k\times n}_+ \, :
\,\, \mathbf{Q}(ilj)=Q_-(il)Q_+(lj),  \big.\\  & \,\, \big.
\qquad\qquad\qquad\qquad  Q_- ,\,\, Q_+ \ge 0,\,\,\, Q_+e=e, \,\,
 e^\top Q_-e=1 \big\}, \\ \\
\mathcal{Q}  = & \left\{ Q\in \mathbb{R}^{m\times n}_+ \, : \,\,
Q(ij)=\sum_l\mathbf{Q}(ilj) \quad {\rm for \,\,  some} \,\,
\mathbf{Q}\in \bq \right\}.
\end{align*}

\medskip
\noindent The interpretation of the sets $\bp, \bq, \mathcal{Q}$
is given next.

\medskip
\noindent Suppose one is given random variables $(Y_-,X,Y_+)$,
taking values in $\{1,\dots ,m\}\times \{1,\dots
,k\}\times\{1,\dots ,n\}$. For convenience we can think of the
r.v.'s as defined on the canonical measurable space
$(\Omega,\mathcal{F})$, where $\Omega$ is the set of all triples
$(i,l,j)$ and $\mathcal{F}$ is $2^\Omega$. For $\omega=(i,l,j)$ we
have the identity mapping $(Y_-,X,Y_+)(\omega)=(i,l,j)$. If
$\mathbb{R}$ a given probability measure on this space, then the
distribution of the triple $(Y_-,X,Y_+)$ under $\mathbb{R}$ is
given by the {\em tensor} $\mathbf{R}$ defined by
\begin{equation}\label{eq:probtens}
\mathbf{R}(ilj)=\mathbb{R}(Y_-=i,X=l,Y_+=j).
\end{equation}
Conversely, a given tensor $\mathbf{R}$ defines a probability
measure $\mathbb{R}$ on $(\Omega,\mathcal{F})$. We will use the
notation $D$ both for I-divergence between tensors and matrices
and for the Kullback-Leibler divergence between probabilities. If
$\mathbf{P}$, $\mathbf{Q}$ are tensors related to probability
measures $\mathbb{P}$ and $\mathbb{Q}$ like in~(\ref{eq:probtens})
we obviously have
$D(\mathbf{P}||\mathbf{Q})=D(\mathbb{P}||\mathbb{Q})$.

\medskip \noindent The sets $\bp, \bq$ correspond to subsets of the
set of all measures on $(\Omega,\mathcal{F})$. In particular $\bp$
corresponds to the subset of all measures whose $Y=(Y_-, Y_+)$
marginal coincides with the given $P$, while $\bq$ corresponds to
the subset of measures under which $Y_-$ and $Y_+$ are
conditionally independent given $X$. The first assertion is
evident by the definition of $\bp$. To prove the second assertion
notice that if
$\mathbb{Q}(Y_-=i,X=l,Y_+=j)=\mathbf{Q}(ilj)=Q_-(il)Q_+(lj)$, then
summing over $j$ one gets $\mathbb{Q}(Y_-=i,X=l)=Q_-(il)$ (since
$Q_+e=e$) and similarly $\mathbb{Q}(Y_+=j|X=l)= Q_+(lj)$. It
follows that $
\mathbb{Q}(Y_-=i,X=l,Y_+=j)=\mathbb{Q}(Y_-=i,X=l)\mathbb{Q}(Y_+=j|X=l)$
 which is equivalent to
\[
\mathbb{Q}(Y_-=i,Y_+=j|X=l)=\mathbb{Q}(Y_-=i|X=l)\mathbb{Q}(Y_+=j|X=l)
\]
i.e. $Y_-, Y_+$ are conditionally independent given $X$.

\medskip
\noindent Finally the set $\mathcal{Q}$ is best interpreted
algebraically as the set of $m\times n$ probability matrices that
admit exact NMF of size $k$.

\medskip
The following observation (taken from~\cite{ps}) motivates our
approach.

\begin{lem}\label{lemma:61}
$P$ admits exact factorization of inner size $k$ iff
$\bp\cap\bq\neq\emptyset$.
\end{lem}
{\bf Proof.} If $\bp\cap\bq\neq\emptyset$ then there exists a
matrix  $\mathbf{Q}\in\bq$ which also belongs to $\bp$, therefore
$P = Q_-Q_+$. Conversely, if we have $P=Q_-Q_+$ with inner size
$k$, then the tensor $\mathbf{P}$ given by
$\mathbf{P}(ilj)=Q_-(il)Q_+(lj)$ clearly belongs to $\bp$. As in
section~\ref{sec:problem} we can w.l.o.g.\ assume that $Q_+e=e$,
so that $\mathbf{P}$ belongs to $\bq$ as well.~\hfill $\square$

\medskip
\noindent We are now ready to give a natural probabilistic
interpretation to the exact NMF problem. The probability matrix
$P$ admits exact NMF $P=Q_-Q_+$ iff there exists at least one
measure on $(\Omega,\mathcal{F})$ whose $Y=(Y_-,Y_+)$ marginal is
$P$ and at the same time making $Y_-$ and $Y_+$ conditionally
independent given $X$.

\bigskip
\noindent Having shown that the exact NMF factorization $P=Q_-Q_+$
is equivalent to $\bp\cap\bq\neq\emptyset$ it is not surprising
that the approximate NMF, corresponding to $\bp\cap\bq
=\emptyset$, can be viewed as a double minimization over the sets
$\bp$ and $\bq$.
\begin{prop}\label{prop:pqq}
Let $P$ be given. The function $(\mathbf{P},\mathbf{Q})\mapsto
D(\mathbf{P}||\mathbf{Q})$ attains a minimum on $\bp\times \bq$
and it holds that
\begin{equation*}
\min_{Q\in\mathcal{Q}}D(P||Q)=\min_{\mathbf{P}\in\bp,\mathbf{Q}\in\bq}D(\mathbf{P}||\mathbf{Q}).
\end{equation*}
\end{prop}

\noindent The proof will be given in
subsection~\ref{subsection:pm}.

\begin{remark}\label{remark:lessk}
{\em Let $\mathbf{P}^*$ and $\mathbf{Q}^*$ be the minimizing
elements in proposition~\ref{prop:pqq}. If there is $l_0$ such
that $\sum_{ij}\mathbf{P}^*(il_0 j)=0$, then all
$\mathbf{Q}^*(il_0j)$ are zero as well. Similarly, if there is
$l_0$ such that $\sum_{ij}\mathbf{Q}^*(i l_0 j)=0$, then all
$\mathbf{P}^*(il_0j)$ are zero as well. In each (and hence both)
of these cases the optimal approximate factorization $Q^*_-Q^*_+$
of $P$ is of inner size less than $k$ (delete the column
corresponding to $l_0$ from $Q^*_-$ and the corresponding row of
$Q^*_+$).}
\end{remark}

\subsection{Two partial minimization
problems}\label{subsection:pm}

In the next section we will construct the algorithm for the
solution of the double minimization problem
$$
\min_{\mathbf{P}\in\bp,\mathbf{Q}\in\bq}D(\mathbf{P}||\mathbf{Q}),
$$
of proposition~\ref{prop:pqq}, as an alternating minimization
algorithm over the two sets $\bp$ and $\bq$. This motivates us to
consider here two partial minimization problems. In the first one,
given $\mathbf{Q}\in\bq$ we minimize the I-divergence
$D(\mathbf{P}||\mathbf{Q})$ over $\mathbf{P}\in\bp$. In the second
problem, given $\mathbf{P} \in\bp$ we minimize the I-divergence
$D(\mathbf{P}||\mathbf{Q})$ over $\mathbf{Q}\in\bq$.

\medskip
Let us start with the first problem. The unique solution
$\mathbf{P}^*=\mathbf{P}^*(\mathbf{Q})$ can easily be computed
analytically and is given by
\begin{equation}\label{eq:p*}
\mathbf{P}^*(ilj)=\frac{\mathbf{Q}(ilj)}{Q(ij)}\, P(ij),
\end{equation}
where $Q(ij)=\sum_l\mathbf{Q}(ilj)$. We also adopt the convention
to put $\mathbf{P}^*(ilj)=0$ if $Q(ij)=0$, which ensures that,
viewed as measures, $\mathbf{P}^*\ll\mathbf{Q}$.

\medskip
Now we turn to the second partial minimization problem. The unique
solution $\mathbf{Q}^*=\mathbf{Q}^*(\mathbf{P})$ to this problem
can also be easily computed analytically and is given by
\begin{align}
Q^*_-(il) & = \sum_j\mathbf{P}(ilj)\label{eq:q-} \\
Q^*_+(lj) & = \frac{\sum_i
\mathbf{P}(ilj)}{\sum_{ij}\mathbf{P}(ilj)},\label{eq:q+}
\end{align}
where we assign arbitrary values to the $Q^*_+(lj)$ (complying
with the constraint $Q_+e=e$) for those $l$ with
$\sum_{ij}\mathbf{P}(ilj)=0$.

\bigskip \noindent The two partial minimization problems and their
solutions have a nice probabilistic interpretation.

\medskip \noindent In the first minimization problem, one is given a distribution
$\mathbf{Q}$, which makes the pair $Y=(Y_-, Y_+)$ conditionally
independent given $X$, and finds the best approximation to it in
the set $\bp$ of distributions with the marginal of $Y$ given by
$P$. Let $\mathbf{P}^*$ denote the optimal distribution of
$(Y_-,X,Y_+)$. Equation~(\ref{eq:p*}) can then be interpreted, in
terms of the corresponding measures, as
$$
\mathbb{P}^*(Y_-=i,X=l,Y_+=j) = \mathbb{Q}(X=l|Y_-=i,Y_+=j)P(ij).
$$
Notice that the conditional distributions of $X$ given $Y$ under
$\mathbb{P}^*$ and $\mathbb{Q}$ are the same. We will see below
that this is not a coincidence.

\medskip
\noindent In the second minimization problem, one is given a
distribution $\mathbf{P}$, with the marginal of $Y$ given by $P$
and finds the best approximation to it in the set $\bq$ of
distributions which make $Y=(Y_-, Y_+)$ conditionally independent
given $X$. Let $\mathbf{Q}^*$ denote the optimal distribution of
$(Y_-,X,Y_+)$. Equations~(\ref{eq:q-}) and~(\ref{eq:q+}) can then
be interpreted, in terms of the corresponding measures, as
\[
\mathbb{Q}^*(Y_-=i,X=l)=\mathbb{P}(Y_-=i,X=l)
\]
and
\[
\mathbb{Q}^*(Y_+=j|X=l)=\mathbb{P}(Y_+=j|X=l).
\]
We see that the optimal solution $\mathbb{Q}^*$ is such that the
marginal distributions of $(X,Y_-)$ under $\mathbb{P}$ and
$\mathbb{Q}^*$ coincide as well as the conditional distributions
of $Y_+$ given $X$ under $\mathbb{P}$ and $\mathbb{Q}^*$. Again,
this is not a coincidence, as we will explain below.

\medskip
\begin{remark}
{\em As a side remark we notice that the minimization of
$D(\mathbf{Q}||\mathbf{P})$ over $\mathbf{P}\in\bp$ for a given
$\mathbf{Q}\in\bq$ yields the same solution $\mathbf{P}^*$. A
similar result does not hold for the second minimization problem.
This remark is not relevant for what follows.}
\end{remark}

\noindent We can now state the so called {\em Pythagorean rules}
for the two partial minimization problems. This terminology was
introduced by Csisz\'ar \cite{c1975}.
\begin{lem}\label{lemma:pyth}
For fixed $\mathbf{Q}$ and $\mathbf{P}^*=\mathbf{P}^*(\mathbf{Q})$
it holds that, for any $\mathbf{P} \in\bp$,
\begin{equation}\label{eq:pythp}
D(\mathbf{P}||\mathbf{Q})=D(\mathbf{P}||\mathbf{P}^*)+D(\mathbf{P}^*||\mathbf{Q}),
\end{equation}
moreover
\begin{equation}\label{eq:p0q0}
D(\mathbf{P}^*||\mathbf{Q})=D(P||Q),
\end{equation}
where
\begin{equation}\label{eq:qq}
Q(ij)=\sum_l\mathbf{Q}(ilj).
\end{equation}
For fixed $\mathbf{P}$ and $\mathbf{Q}^*=\mathbf{Q}^*(\mathbf{P})$
it holds that, for any $\mathbf{Q} \in\bq$,
\begin{equation}\label{eq:pythq}
D(\mathbf{P}||\mathbf{Q})=D(\mathbf{P}||\mathbf{Q}^*)+D(\mathbf{Q}^*||\mathbf{Q}).
\end{equation}
\end{lem}
{\bf Proof.} To prove the first rule we compute
\begin{eqnarray*}
\lefteqn{D(\mathbf{P}||\mathbf{P}^*)+ D(\mathbf{P}^*||\mathbf{Q})}\\
& = &
\sum_{ilj}\mathbf{P}(ilj)\log\frac{\mathbf{P}(ilj)Q(ij)}{\mathbf{Q}(ilj)P(ij)}
+
\sum_{ilj}\mathbf{Q}(ilj)\frac{P(ij)}{Q(ij)}\log \frac{P(ij)}{Q(ij)} \\
& = &
\sum_{ilj}\mathbf{P}(ilj)\log\frac{\mathbf{P}(ilj)}{\mathbf{Q}(ilj)}
+
\sum_{ilj}\mathbf{P}(ilj)\log\frac{Q(ij)}{P(ij)} \\
& & \mbox{} + \sum_{ij}Q(ij)\frac{P(ij)}{Q(ij)}\log
\frac{P(ij)}{Q(ij)} = D(\mathbf{P}||\mathbf{Q}).
\end{eqnarray*}
The first rule follows. To prove the relation~(\ref{eq:p0q0})
insert equation~(\ref{eq:p*}) into $D(\mathbf{P}^*||\mathbf{Q})$
and sum over $l$ to get
\[
D(\mathbf{P}^*||\mathbf{Q})=\sum_{ilj}P(ij)\frac{\mathbf{Q}(ilj)}{Q(ij)}\log
\frac{P(ij)}{Q(ij)} = D(P||Q).
\]

\noindent To prove the second rule we first introduce some
notation. Let $\mathbf{P}(il\cdot)=\sum_j \mathbf{P}(ilj)$,
$\mathbf{P}(\cdot lj)=\sum_i\mathbf{P}(ilj)$ and
$\mathbf{P}(j|l)=\mathbf{P}(\cdot lj)/\sum_j\mathbf{P}(\cdot lj)$.
For $\mathbf{Q}$ we use similar notation and observe that
$\mathbf{Q}(il\cdot)=Q_-(il)$, and
$\mathbf{Q}(j|l)=Q_+(lj)/\sum_jQ_+(lj)$, and
$Q^*_-(il)=\mathbf{P}(il\cdot)$ and $Q^*_+(lj)=\mathbf{P}(j|l)$.
We now compute
\begin{align*}
D(\mathbf{P}||\mathbf{Q})-D(\mathbf{P}||\mathbf{Q}^*)
&=\sum_{ilj}\mathbf{P}(ilj) \left( \log
\frac{\mathbf{P}(il\cdot)}{Q_-(il)}+\log\frac{\mathbf{P}(j|l)}{Q_+(lj)} \right)\\
& = \sum_{il}\mathbf{P}(il\cdot)\log
\frac{\mathbf{P}(il\cdot)}{Q_-(il)} + \sum_{lj}\mathbf{P}(\cdot
lj)\log\frac{\mathbf{P}(j|l)}{Q_+(lj)} \\
& = D(\mathbf{Q}^*||\mathbf{Q}).
\end{align*}
The second rule follows. \hfill$\square$

\bigskip \noindent With the aid of the relation~(\ref{eq:p0q0}) we can now
prove
proposition~\ref{prop:pqq}. \medskip\\
{\bf Proof of proposition~\ref{prop:pqq}.} With
$\mathbf{P}^*=\mathbf{P}^*(\mathbf{Q})$, the optimal solution of
the partial minimization over $\bp$, we have
\begin{align*}
D(\mathbf{P}||\mathbf{Q})& \geq D(\mathbf{P}^*||\mathbf{Q}) \\
& = D(P||Q) \\
& \geq  \min_{Q\in\mathcal{Q}}D(P||Q).
\end{align*}
It follows that
$\inf_{\mathbf{P}\in\bp,\mathbf{Q}\in\bq}D(\mathbf{P}||\mathbf{Q})\geq\min_{Q\in\mathcal{Q}}D(P||Q)$.
\\
Conversely, let $\mathbf{Q}$  in $\bq$ be given and let $Q$ be
defined by $Q(ij)=\sum_l\mathbf{Q}(ilj)$ . From
\begin{align*}
D(P||Q) & = D(\mathbf{P}^*(\mathbf{Q})|| \mathbf{Q})\\
& \geq
\inf_{\mathbf{P}\in\bp,\mathbf{Q}\in\bq}D(\mathbf{P}||\mathbf{Q}),
\end{align*}
we obtain
\[
\min_{Q\in\mathcal{Q}}D(P||Q)\geq\inf_{\mathbf{P}\in\bp,\mathbf{Q}\in\bq}D(\mathbf{P}||\mathbf{Q}).
\]
Finally we show that we can replace the infima by minima. Let
$Q^*_-$ and $Q^*_+$ be such that $(Q_-,Q^+)\mapsto D(P||Q_-Q^+)$
is minimized (their existence is guaranteed by
proposition~\ref{prop:exist}). Let $\mathbf{Q}^*$ be a
corresponding element in $\bq$ and
$\mathbf{P}^*=\mathbf{P}^*(\mathbf{Q}^*)$. Then
$D(\mathbf{P}^*||\mathbf{Q}^*)=D(P||Q^*_-Q^*_+)$ and the result
follows.
\hfill$\square$

\bigskip
For a probabilistic derivation of the solutions of the two partial
minimization problems and of their corresponding Py\-tha\-go\-rean
rules, we use a general result (lemma~\ref{lemma:uv} below) on the
I-divergence between two joint laws of any random vector $(U,V)$.
We denote the law of $(U,V)$ under arbitrary probability measures
$\mathbb{P}$ and $\mathbb{Q}$ by $\mathbb{P}^{U,V}$ and
$\mathbb{Q}^{U,V}$. The conditional distributions of $U$ given $V$
are summarized by the matrices $\mathbb{P}^{U|V}$ and
$\mathbb{Q}^{U|V}$, with the obvious convention
$\mathbb{P}^{U|V}(ij)=\mathbb{P}(U=j|V=i)$ and likewise for
$\mathbb{Q}^{U|V}$.
\begin{lem}\label{lemma:uv}
It
holds that
\begin{equation}\label{eq:duv}
D(\mathbb{P}^{U,V}||\mathbb{Q}^{U,V})=\mathbb{E}_\mathbb{P}
D(\mathbb{P}^{U|V}||\mathbb{Q}^{U|V}) +
D(\mathbb{P}^V||\mathbb{Q}^V),
\end{equation}
where
\[
D(\mathbb{P}^{U|V}||\mathbb{Q}^{U|V}) = \sum_j
P(U=j|V)\log\frac{P(U=j|V)}{Q(U=j|V)}.
\]
If moreover $V=(V_1,V_2)$, and $U, V_2$ are conditionally
independent given $V_1$ under $\mathbb{Q}$, then the first term on
the RHS of~(\ref{eq:duv}) can be written as
\begin{equation}\label{eq:duv1} \mathbb{E}_\mathbb{P}
D(\mathbb{P}^{U|V}||\mathbb{Q}^{U|V}) =\mathbb{E}_\mathbb{P}
D(\mathbb{P}^{U|V}||\mathbb{P}^{U|V_1}) + \mathbb{E}_\mathbb{P}
D(\mathbb{P}^{U|V_1}||\mathbb{Q}^{U|V_1}).
\end{equation}
\end{lem}
{\bf Proof.} It follows from elementary manipulations.
\hfill $\square$\medskip\\
The first minimization problem can be solved probabilistically as
follows. Given $\mathbf{Q}$ we are to find its best
approximation within $\bp$. Let $\mathbb{Q}$ correspond to the
given $\mathbf{Q}$ and $\mathbb{P}$ correspond to the generic
$\mathbf{P} \in \bp$. Choosing $U=X$, $V=Y=(Y_-,Y_+)$ in lemma
\ref{lemma:uv}, and remembering that $\mathbb{P}^Y$ is determined
by $P$ for all $\mathbf{P} \in \bp$, equation (\ref{eq:duv}) now
reads
\begin{equation}\label{eq:pq1}
D(\mathbf{P}||\mathbf{Q})=\mathbb{E}_{\mathbb{P}}
D(\mathbb{P}^{X|Y}||\mathbb{Q}^{X|Y}) + D(P||Q),
\end{equation}
where the matrix $Q$ is as in~(\ref{eq:qq}).
The problem is equivalent to the minimization
of $\mathbb{E}_{\mathbb{P}}D(\mathbb{P}^{X|Y}||\mathbb{Q}^{X|Y})$
w.r.t. $\mathbf{P} \in \bp$, which is attained (with value $0$) at
$\mathbb{P}^*$ with $\mathbb{P}^{*\, X|Y}=\mathbb{Q}^{X|Y}$ and
$\mathbb{P}^{*Y}= P$. To derive probabilistically the
corresponding Pythagorean rule, we apply~(\ref{eq:duv}) with
$\mathbb{P}^*$ instead of $\mathbb{Q}$. We obtain, using
$\mathbb{P}^Y=\mathbb{P}^{*Y}$,
\begin{equation}\label{eq:pq8}
D(\mathbb{P}^{X,Y}||\mathbb{P}^{*^{X,Y}})=
\mathbb{E}_\mathbb{P}D(\mathbb{P}^{X|Y}||\mathbb{P}^{*^{X|Y}}).
\end{equation}
Since also
\begin{equation}\label{eq:pq2}
\mathbb{E}_\mathbb{P}D(\mathbb{P}^{X|Y}||\mathbb{Q}^{X|Y})=\mathbb{E}_\mathbb{P}
D(\mathbb{P}^{X|Y}||\mathbb{P}^{*^{X|Y}}),
\end{equation}
we combine equations~(\ref{eq:pq8}) and~(\ref{eq:pq2}) and insert
the result into~(\ref{eq:pq1}). Recognizing the fact that
$D(\mathbf{P}||\mathbf{P}^*)=D(\mathbb{P}^{X,Y}||\mathbb{P}^{*^{X,Y}})$,
and using $D(\mathbf{P}^*||\mathbf{Q})=D(P||Q)$
according to (\ref{eq:p0q0}), we then identify~(\ref{eq:pq1}) as
the first Pythagorean rule~(\ref{eq:pythp}).
\medskip\\
The treatment of the second minimization problem follows a similar
pattern. Given $\mathbf{P}$ we are to find its best
approximation within $\bq$.  Let $\mathbb{P}$ correspond to the
given $\mathbf{P}$ and $\mathbb{Q}$ correspond to the generic
$\mathbf{Q} \in \bq$. Choosing $U=Y_+$, $V_1=X$ and $V_2=Y_-$ in
lemma \ref{lemma:uv}, and remembering that under any
$\mathbf{Q}\in\bq$ the r.v. $Y_-, Y_+$ are conditionally
independent given $X$, equation (\ref{eq:duv}) refined with
(\ref{eq:duv1}) now reads
\begin{align*}
D(\mathbf{P}||\mathbf{Q})= &
\mathbb{E}_\mathbb{P}D(\mathbb{P}^{Y_+|X,Y_-}||\mathbb{P}^{Y_+|X})\\
& \mbox{}
+\mathbb{E}_\mathbb{P}D(\mathbb{P}^{Y_+|X}||\mathbb{Q}^{Y_+|X})+D(\mathbb{P}^{Y_-,X}||\mathbb{Q}^{Y_-,X}).
\end{align*}
The problem is equivalent to the minimizations of the second and
third I-divergences on the RHS w.r.t. $\mathbf{Q} \in \bq$, which
are attained (both with value $0$) at $\mathbb{Q}^*$ with
$\mathbb{Q}^{*\, Y_+|X}=\mathbb{P}^{Y_+|X}$ and
$\mathbb{Q}^{*Y_-,X}= \mathbb{P}^{Y_-,X}$. Note that $X$ has the
same distribution under $\mathbb{P}$ and $\mathbb{Q}^*$. To derive
probabilistically the corresponding Pythagorean rule we notice
that
\begin{equation}\label{eq:pq3}
D(\mathbf{P}||\mathbf{Q})-D(\mathbf{P}||\mathbf{Q}^*) =
\mathbb{E}_{\mathbb{Q}^*}D(\mathbb{Q}^*{^{Y_+|X}}||\mathbb{Q}^{Y_+|X})+D(\mathbb{Q}^{*^{Y_-,X}}||\mathbb{Q}^{Y_-,X}).
\end{equation}
In the right hand side of~(\ref{eq:pq3}) we can, by conditional
independence, replace
$\mathbb{E}_{\mathbb{Q}^*}D(\mathbb{Q}^*{^{Y_+|X}}||\mathbb{Q}^{Y_+|X})$
with
$\mathbb{E}_{\mathbb{Q}^*}D(\mathbb{Q}^*{^{Y_+|X,Y_-}}||\mathbb{Q}^{Y_+|X,Y-})$.
By yet another application of~(\ref{eq:duv}), we thus see that
$D(\mathbf{P}||\mathbf{Q})-D(\mathbf{P}||\mathbf{Q}^*)=D(\mathbf{Q}^*||\mathbf{Q})$,
which is the second Pythagorean rule~(\ref{eq:pythq}).

\section{Alternating minimization algorithm}\label{section:altmin}

The results of the previous section are aimed at setting up an
alternating minimization algorithm for obtaining $\min_Q D(P||Q)$,
where $P$ is a given nonnegative matrix. In view of
proposition~\ref{prop:pqq} we can lift this problem to the
$\bp\times \bq$ space. Starting with an arbitrary $\mathbf{Q}^0
\in\bq$ with positive elements, we adopt the following alternating
minimization scheme
\begin{equation}\label{eq:altalgo}
 \to \mathbf{Q}^t \to \mathbf{P}^t\to \mathbf{Q}^{t+1}\to \mathbf{P}^{t+1}
\to
\end{equation}
where $\mathbf{P}^t=\mathbf{P}^*(\mathbf{Q}^t)$,
$\mathbf{Q}^{t+1}=\mathbf{Q}^*(\mathbf{P}^t)$.
\medskip\\
To relate this algorithm to the one of section~\ref{sec:problem}
(formulas (\ref{eq:q-bar2}) and (\ref{eq:q+bar2})) we combine two
steps of the alternating minimization at a time.
From~(\ref{eq:altalgo}) we get
$$\mathbf{Q}^{t+1}=\mathbf{Q}^*(\mathbf{P}^*(\mathbf{Q}^t)).$$
\noindent Computing the optimal solutions according
to~(\ref{eq:p*}),~(\ref{eq:q-}) and (\ref{eq:q+}) one gets from
here the formulas (\ref{eq:q-bar2}) and (\ref{eq:q+bar2}) of
section~\ref{sec:problem}.

\medskip
\noindent The Pythagorean rules allow us to easily compute the
update gain $D(P||Q^t)-D(P||Q^{t+1})$ of the algorithm.
\begin{prop}\label{prop:gain}
The update gain at each iteration of the
algorithm~(\ref{eq:altalgo}) in terms of the matrices $Q^t$ is
given by
\begin{equation}\label{eq:gain}
D(P||Q^{t}) - D(P||Q^{t+1}) =
D(\mathbf{P}^{t}||\mathbf{P}^{t+1})+D(\mathbf{Q}^{t+1}||\mathbf{Q}^t).
\end{equation}
\end{prop}
{\bf Proof.} The two Pythagorean rules from lemma~\ref{lemma:pyth}
now take the forms
\begin{align*}
D(\mathbf{P}^t||\mathbf{Q}^t) & =
D(\mathbf{P}^t||\mathbf{Q}^{t+1})+D(\mathbf{Q}^{t+1}||\mathbf{Q}^t),
\\
D(\mathbf{P}^{t}||\mathbf{Q}^{t+1}) & =
D(\mathbf{P}^{t}||\mathbf{P}^{t+1})+D(\mathbf{P}^{t+1}||\mathbf{Q}^{t+1}).
\end{align*}
Addition  of these two equations
results in
\begin{align*}
D(\mathbf{P}^{t}||\mathbf{Q}^t) & =
D(\mathbf{P}^{t}||\mathbf{P}^{t+1})+D(\mathbf{P}^{t+1}||\mathbf{Q}^{t+1})+D(\mathbf{Q}^{t+1}||\mathbf{Q}^t),
\end{align*}
and since $D(\mathbf{P}^{t}||\mathbf{Q}^t) = D(P||Q^{t})$ from
(\ref{eq:p0q0}), the result follows. \hfill$\square$

\begin{remark}
{\em If one starts the algorithm with matrices $(Q^0_-, Q^0_+)$ in
the interior of the domain, the iterations will remain in the
interior. Suppose that, at step $n$, the update gain is zero.
Then, from~(\ref{eq:gain}), we get that
$D(\mathbf{Q}^{t+1}||\mathbf{Q}^t)=0$. Hence the tensors
$\mathbf{Q}^{t+1}$ and $\mathbf{Q}^t$ are identical. From this it
follows by summation that $Q^{t+1}_-=Q^t_-$. But then we also have
the equality $Q^t_-(il)Q^{t+1}_+(lj)=Q^t_-(il)Q^{t}_+(lj)$ for all
$i,l,j$. Since all $Q^t_-(il)$ are positive, we also have
$Q^{t+1}_+=Q^t_+$. Hence, the updating formulas strictly decrease
the objective function until the algorithm reaches a fixed point.}
\end{remark}
\rm
We close this section with the proof of
proposition~\ref{prop:exist} in which we use the result of
proposition~\ref{prop:gain}.
\medskip\\
{\bf Proof of proposition~\ref{prop:exist}.} We first prove that
there exists a pair of matrices $(W,H)$ with $He_m=e_k$ and
$We_k=Ve_n$ for which $D(V||WH)$ is finite. Put
$W=\frac{1}{k}Ve_ne_k^\top$ and $H=\frac{1}{e_m^\top Ve_n}e_k
e_m^\top V$. Note that indeed $He_m=e_k$ and $We_k=Ve_n$ and that
all elements of $W$ and $H$, and hence those of $WH$, are
positive, $D(V||WH)$ is therefore finite.

Next we show that we can restrict ourselves to minimization over a
compact set $\mathcal{K}$ of matrices. Specifically, we will show
that for all positive matrices $W$ and $H$, there exist positive
matrices $W'$ and $H'$ with $(W',H')\in \mathcal{K}$ such that
$D(V||W'H')\leq D(V||WH)$. We choose for arbitrary $W^0$ and $H^0$
the matrices $W^1$ and $H^1$ according to~(\ref{eq:wbar2})
and~(\ref{eq:hbar2}). It follows from proposition~\ref{prop:gain}
that indeed $D(V||W^1H^1)\leq D(V||W^0H^0)$. Moreover, it is
immediately clear from~(\ref{eq:wbar2}) and~(\ref{eq:hbar2}) that
we have $W^1e=Ve$ and $H^1e=e$. Hence, it is sufficient to confine
search to the compact set $\mathcal{L}$ where $He=e$ and $We=Ve$.

Fix a pair of indices $i,j$. Since we can compute the divergence
elementwise we have the trivial estimate
\[
D(V||WH)\geq V_{ij}\log \frac{V_{ij}}{(WH)_{ij}}-V_{ij}+(WH)_{ij}.
\]
Since for $V_{ij}>0$ the function $d_{ij}:x\to V_{ij}\log
\frac{V_{ij}}{x}-V_{ij}+x$ is  decreasing on $(0,V_{ij})$, we have
for any sufficiently small $\eps>0$ (of course $\eps < V_{ij}$)
that $d_{ij}(x)>d_{ij}(\eps)$ for $x\leq \eps$ and of course
$\lim_{\eps\to 0}d_{ij}(\eps)=\infty$. Hence to find the minimum
of $d_{ij}$, it is sufficient to look at $x\geq \eps$. Let $\eps_0
>0$ and such that $\eps_0 <\min\{V_{ij}:V_{ij}>0\}$. Let $\mathcal{G}$ be
the set of $(W,H)$ such that $(WH)_{ij}\geq \eps_0$ for all $i,j$
with $V_{ij}>0$. Then $\mathcal{G}$ is closed. Take now
$\mathcal{K}=\mathcal{L}\cap \mathcal{G}$, then $\mathcal{K}$ is
the compact set we are after. Let us observe that $\mathcal{K}$ is
non-void for sufficiently small $\eps_0$. Clearly the map
$(W,H)\mapsto D(V||WH)$ is continuous on $\mathcal{K}$ and thus
attains its minimum. \hfill\ $\square$

\section{Auxiliary functions}\label{section:auxfunc}

Algorithms for recursive minimization can often be constructed by
using {\em auxiliary functions}. For the problem of minimizing the
divergence $D(V||WH)$, some such functions can be found
in~\cite{leeseung2001} and they are analogous to functions that
are used when studying the EM algorithm, see~\cite{wu}. The choice
of an auxiliary function is usually based on {\em ad hoc}
reasoning, like for instance finding a Lyapunov function for
studying the stability of the solutions of a differential
equation. We show in this section that the lifted version of the
divergence minimization problem leads in a natural way to useful
auxiliary functions. Let us first explain what is meant by an
auxiliary function.

Suppose one wants to minimize a function $x\mapsto F(x)$, defined
on some domain. The function $(x,x')\mapsto G(x,x')$ is an
auxiliary function for $F$ if
\begin{align*}
G(x,x')  & \geq  F(x'), \,\,\, \forall x,x',\\
G(x,x) & =  F(x), \quad \forall x.
\end{align*}
If we define (assuming that the $\arg\min$ below exists and is
unique)
\begin{equation}\label{eq:update}
x'=x'(x)=\arg\min G(x,\cdot),
\end{equation}
then we have
\[
F(x')\leq G(x,x') \leq G(x,x)=F(x),
\]
and hence the value of $F$ decreases by replacing $x$ with $x'$. A
recursive procedure to find the minimum of $F$ can be based on the
recipe~(\ref{eq:update}) by taking $x=x^t$ and $x'=x^{t+1}$. To be
useful an auxiliary function $G$ must allow for a simple
computation or characterization of $\arg\min G(x,\cdot)$.

\medskip
We consider now the minimization of $D(P||Q)$ and its lifted
version, the minimization of $D(\mathbf{P}||\mathbf{Q})$ as in
section~\ref{section:lift}. In particular, with reference to the
alternating minimization scheme~(\ref{eq:altalgo}), with the
notations of section~\ref{section:altmin}, we know that
$\mathbf{Q}^{t+1}$ is found by minimizing $\mathbf{Q}'\mapsto
D(\mathbf{P}^*(\mathbf{Q}^t)||\mathbf{Q}')$. This strongly
motivates the choice of the function
\[ (\mathbf{Q},\mathbf{Q}') \mapsto G(\mathbf{Q},\mathbf{Q}')=D(\mathbf{P}^*(\mathbf{Q})||\mathbf{Q}')
\]
as an auxiliary function for minimizing $D(P||Q)$ w.r.t. $Q$.

\medskip
Using the decomposition of the divergence in
equation~(\ref{eq:duv}) we can rewrite $G$ as
\begin{equation}\label{eq:div29}
G(\mathbf{Q},\mathbf{Q}')=D(\mathbb{P}^{*^Y}||\mathbb{Q}'^Y)+\mathbb{E}_{\mathbb{P}^*}D(\mathbb{P}^{*^{X|Y}}||\mathbb{Q}'^{X|Y}).
\end{equation}
Since $\mathbb{P}^{*X|Y} = \mathbb{Q}^{X|Y}$, and $\mathbb{P}^{*Y}
= P$ we can rewrite~(\ref{eq:div29}) as
\begin{equation}\label{eq:div30}
G(\mathbf{Q},\mathbf{Q}')=D(P||\mathbb{Q}'^Y)+\mathbb{E}_PD(\mathbb{Q}^{X|Y}||\mathbb{Q}'^{X|Y}).
\end{equation}
From~(\ref{eq:div30}) it follows that
$G(\mathbf{Q},\mathbf{Q}')\geq D(P||Q')$, and that
$G(\mathbf{Q},\mathbf{Q})= D(P||Q)$, precisely the two properties
that define an auxiliary function for $D(P||Q)$.

\medskip\noindent
In~\cite{leeseung2001} one can find two auxiliary functions for
the original minimization problem $D(V||WH)$. One function is for
minimization over $H$ with fixed $W$, the other for minimization
over $W$ with fixed $H$. To show the connection with the function
$G$ defined above, we first make the dependence of $G$ on
$Q_-,Q_+,Q'_-,Q'_+$ explicit by writing
$G(\mathbf{Q},\mathbf{Q}')$ as $G(Q_-,Q_+,Q'_-,Q'_+)$.

\medskip \noindent The auxiliary function for minimization with fixed
$Q_-$ can then be taken as \[Q'_+\mapsto
G^+_{\mathbf{Q}}(Q'_+)=G(Q_-,Q_+,Q_-,Q'_+),\]

\noindent whereas the auxiliary function for minimization with
fixed $Q_+$ can be taken as \[Q'_-\mapsto
G^-_{\mathbf{Q}}(Q'_-)=G(Q_-,Q_+,Q'_-,Q_+)\]

\medskip \noindent The functions $G^+_{\mathbf{Q}}$ and $G^-_{\mathbf{Q}}$
correspond to the auxiliary functions in~\cite{leeseung2001},
where they are given in an explicit form, but where no rationale
for them is given.

For the different auxiliary functions introduced above, we will
now compute the update gains and compare these expressions with
(\ref{eq:gain}).
\begin{lem}\label{lemma:auxdiff}
Consider the auxiliary functions $G$ ,$G^-_{\mathbf{Q}}$ ,
$G^+_{\mathbf{Q}}$ above. Denote by $Q'_-$ and $Q'_+$ the
minimizers of the auxiliary functions in all three cases. The
following equalities hold
\begin{align}
D(P||Q_-Q_+)-G^-_{\mathbf{Q}}(Q'_-)  & =
D(\mathbb{Q}'^{Y_-,X}||\mathbb{Q}^{Y_-,X}) \label{eq:auxdiff1}\\
D(P||Q_-Q_+)-G^+_{\mathbf{Q}}(Q'_+)  & =
\mathbb{E}_{\mathbb{P}^*}D(\mathbb{Q}'^{Y_+|X}||\mathbb{Q}^{Y_+|X})\label{eq:auxdiff2}\\
D(P||Q_-Q_+)-G(Q_-,Q_+,Q'_-,Q'_+) & =
D(\mathbb{Q}'^{Y_-,X}||\mathbb{Q}^{Y_-,X}) \nonumber \\
&  \,\, \mbox{~~~}
+\mathbb{E}_{\mathbb{Q}'}D(\mathbb{Q}'^{Y_+|X}||\mathbb{Q}^{Y_+|X}).
\label{eq:auxdiff3}
\end{align}
\end{lem}
{\bf Proof.} We prove~(\ref{eq:auxdiff3}) first. The other two
follow from this. A simple computation, valid for any $Q_-$a nd $Q_+$, yields
\begin{align}
\lefteqn{D(P||Q_-Q_+)-G(Q_-,Q_+,Q'_-,Q'_+)} \\
&
=\sum_{ij}P(ij)\sum_l\frac{\mathbf{Q}(ilj)}{Q(ij)}\left(\log\frac{Q'_-(il)}{Q_-(il)}
+\log\frac{Q'_+(lj)}{Q_+(lj)}\right) \nonumber\\
& =
\sum_{il}\big(\sum_j\frac{P(ij)\mathbf{Q}(ilj)}{Q(ij)}\big)\log\frac{Q'_-(il)}{Q_-(il)}
 +
\sum_{lj}\big(\sum_i\frac{P(ij)\mathbf{Q}(ilj)}{Q(ij)}\big)\log\frac{Q'_+(lj)}{Q_+(lj)}\label{eq:dg}
\end{align}
Now we exploit the known formulas~(\ref{eq:q-bar2})
and~(\ref{eq:q+bar2}) for the optimizing $Q'_-$ and $Q'_+$. The
first term in~(\ref{eq:dg}) becomes in view of~(\ref{eq:q-bar2})
(or, equivalently, in view of~(\ref{eq:p*}) and~(\ref{eq:q-}))
\[
\sum_{il}Q'_-(il)\log\frac{Q'_-(il)}{Q_-(il)},
\]
which gives the first term on the RHS of~(\ref{eq:auxdiff3}). Similarly, the
second term in~(\ref{eq:dg}) can be written in view
of~(\ref{eq:q+bar2}) as
\[
\sum_{l}\big(\sum_{ij}\mathbf{Q}'(ilj)\big)\sum_{j}Q'_+(lj)\log\frac{Q'_+(lj)}{Q_+(lj)},
\]
which yields the second term on the RHS of
formula~(\ref{eq:auxdiff3}). Formulas~(\ref{eq:auxdiff1})
and~(\ref{eq:auxdiff2}) are obtained similarly, noticing that
optimization of $G^+_{\mathbf{Q}}$ and $G^-_{\mathbf{Q}}$
separately yield the same $Q'_+$, respectively $Q'_-$, as those
obtained by minimization of $G$.
 \hfill$\square$
\begin{remark}
{\em Notice that although for instance $G^-_{\mathbf{Q}}(Q'_-)\geq
D(P||Q'_-Q'_+)$ for all $Q'_-$ and $Q'_+$, we have for the optimal
$Q'_-$ that $G^-_{\mathbf{Q}}(Q'_-)\leq D(P||Q_-Q_+)$.
}
\end{remark}
\begin{cor}\label{cor:corgain}
The update gain of the
algorithm~(\ref{eq:q-bar2}),~(\ref{eq:q+bar2}) can be represented
by
\begin{align}\label{eq:corgain}
D(P||Q^t)-D(P||Q^{t+1})  =
& D(\mathbb{Q}^{{t+1}^{Y_-,X}}||\mathbb{Q}^{t^{Y_-,X}}) \nonumber\\
& +
\mathbb{E}_{\mathbb{Q}^{t+1}}D(\mathbb{Q}^{{t+1}^{Y_+|X}}||\mathbb{Q}^{t^{Y_+|X}})
\nonumber\\
& + \mathbb{E}_PD(\mathbb{Q}^{t^{X|Y}}||\mathbb{Q}^{{t+1}^{X|Y}}).
\end{align}
\end{cor}
{\bf Proof.} Write
\begin{align*}
\lefteqn{D(P||Q^t)-D(P||Q^{t+1})=} \\ & \mbox{~~~~}
D(P||Q^t)-G(\mathbf{Q}^t,\mathbf{Q}^{t+1})+G(\mathbf{Q}^t,\mathbf{Q}^{t+1})-D(P||Q^{t+1})
\end{align*}
and use equations~(\ref{eq:div29}) and~(\ref{eq:auxdiff3}).
\hfill$\square$
\medskip\\
We return to the update formula~(\ref{eq:gain}). A computation
shows the following equalities.
\begin{align}
D(\mathbf{P}^t||\mathbf{P}^{t+1})= &
\mathbb{E}_{P}D(\mathbb{Q}^{t^{X|Y}}||\mathbb{Q}^{{t+1}^{X|Y}})\label{eq:gain1}
\\
D(\mathbf{Q}^{t+1}||\mathbf{Q}^t) = &
D(\mathbb{Q}^{{t+1}^{Y_-,X}}||\mathbb{Q}^{{t}^{Y_-,X}}) \nonumber\\
& +
\mathbb{E}_{\mathbb{Q}^{t+1}}D(\mathbb{Q}^{{t+1}^{Y_+|X}}||\mathbb{Q}^{{t}^{Y_+|X}}).\label{eq:gain2}
\end{align}
In equation~(\ref{eq:gain1}) we recognize the second term in the
auxiliary function, see~(\ref{eq:div30}).
Equation~(\ref{eq:gain2}) corresponds to
equation~(\ref{eq:auxdiff3}) of lemma~\ref{lemma:auxdiff} and
we see that formula~(\ref{eq:gain}) is indeed the same
as~(\ref{eq:corgain}) .
\medskip\\
The algorithm~(\ref{eq:q-bar2}),~(\ref{eq:q+bar2}) is to be
understood by using these two equations simultaneously. As an
alternative one could first use~(\ref{eq:q-bar2}) to obtain
$Q^{t+1}_-$ and, instead of using $Q^t_-$, feed this result
into~(\ref{eq:q+bar2}) to obtain $Q^{t+1}_+$. If we do this, we
can express the update gain of the first partial step, like in the
proof of corollary~\ref{cor:corgain}, by adding the result of
equation~(\ref{eq:auxdiff1}) to the second summand
of~(\ref{eq:div30}), with the understanding that $\mathbb{Q}'$ is
now given by the $Q^{t+1}(ij)Q^t(lj)$. The update gain of the
second partial step is likewise obtained by combining the result
of~(\ref{eq:auxdiff2}) and the second summand of~(\ref{eq:div30}),
with the understanding that now $\mathbb{Q}$ is to be interpreted
as given by the $Q^{t+1}(ij)Q^t(lj)$. Of course, as another
alternative, the order of the partial steps can be reversed.
Clearly, the expressions for the update gains for these cases also
result from working with the auxiliary functions
$G^-_{\mathbf{Q}}$ and $G^+_{\mathbf{Q}}$, the
equations~(\ref{eq:auxdiff1}) and~(\ref{eq:auxdiff2}) and
proceeding as in the proof of corollary~\ref{cor:corgain}.

\section{Convergence properties}\label{section:convprop}

In this section we study the convergence properties of the
divergence minimization algorithm (\ref{eq:q-bar2}),
(\ref{eq:q+bar2}).

The  next theorem states that the sequences generated by the
algorithm converge for every (admissible) initial value. Of course
the limits will in general depend on the initial value.

\begin{thm}\label{thm:conv}
Let $Q^t_-(il)$, $Q_+^t(lj)$ be generated by the algorithm
(\ref{eq:q-bar2}), (\ref{eq:q+bar2}) and $\mathbf{Q}^t$ the
corresponding tensors. Then the $Q^t_-(il)$ converge to  limits
$Q_-^\infty(il)$ and the $\mathbf{Q}^t$ converges to a limit
$\mathbf{Q}^\infty$ in $\bq$. The $Q_+^t(lj)$ converge to limits
$Q_+^\infty(lj)$ for all $l$ with $\sum_iQ_+^\infty(il)>0$.
\end{thm}
{\bf Proof.} We first show that the $Q^t_-$ and $Q^t_+$ form
convergent sequences. We start with equation~(\ref{eq:gain}). By
summing over $n$ we obtain
\[
D(P||Q^0)-D(P||Q^t)=\sum_{k=1}^{t-1}\Big(D(\mathbf{P}^{s}||\mathbf{P}^{s+1})
+ D(\mathbf{Q}^{s+1}||\mathbf{Q}^s)\Big).
\]
It follows that
$\sum_{k=1}^{\infty}D(\mathbf{P}^{s}||\mathbf{P}^{s+1})$ and
$\sum_{k=1}^{\infty}D(\mathbf{Q}^{s+1}||\mathbf{Q}^s)$ are finite.
Now we use that fact that for any two probability measures, the
Kullback-Leibler divergence $D(\mathbb{P}||\mathbb{Q})$ is greater
than or equal to their Hellinger distance
$H(\mathbb{P},\mathbb{Q})$, which is the $L^2$ distance between
the square roots of corresponding densities w.r.t.~some dominating
measure, see~\cite[p.~368]{shiryaev}. In our case we have
$H(\mathbb{Q}^{s},\mathbb{Q}^{s+1})=\sum_{ilj}
(\sqrt{\mathbf{Q}^{s+1}(ilj)}-\sqrt{\mathbf{Q}^{s}(ilj)})^2$. So
we obtain that
\[
\sum_{k=1}^{\infty}H(\mathbf{Q}^{s+1},\mathbf{Q}^s)<\infty.
\]
We therefore have that, pointwise, the tensors $\mathbf{Q}^t$ form
a Cauchy sequence and hence have a limit $\mathbf{Q}^\infty$. We
will show that $\mathbf{Q}^\infty$ belongs to $\bq$. Since the
$\mathbf{Q}^t(ilj)$ converge to limits $\mathbf{Q}^\infty(ilj)$,
by summation we have that the  marginals
$Q^t_-(il)=\mathbf{Q}^t(il\cdot)$ converge to limits
$\mathbf{Q}^\infty(il\cdot)$ (we use the notation of the proof of
lemma~\ref{lemma:pyth}), and likewise we have convergence of the
marginals $\mathbf{Q}^t(\cdot lj)$ to $\mathbf{Q}^\infty(\cdot
lj)$ and $\mathbf{Q}^t(\cdot l\cdot)$ to $\mathbf{Q}^\infty(\cdot
l\cdot)$. Hence, if $\mathbf{Q}^\infty(\cdot l\cdot)>0$, then  the
$Q^t_+(lj)$ converge to $Q^\infty_+(ij):=\mathbf{Q}^\infty(\cdot
lj)/\mathbf{Q}^\infty(\cdot l\cdot)$ and we have
$\mathbf{Q}^\infty(ilj)=\mathbf{Q}^\infty(il\cdot)Q^\infty_+(ij)$.
Now we analyze the case where $\mathbf{Q}^\infty(\cdot
l_0\cdot)=0$ for some $l_0$. Since in this case both
$\mathbf{Q}^\infty(il_0j)$ and $\mathbf{Q}^\infty(il_0\,\cdot)$
are zero, we have still have a factorization
$\mathbf{Q}^\infty(il_0j)=Q^\infty_-(il_0)Q^\infty_+(l_0j)$, where
we can assign to the $Q^\infty_+(l_0j)$ arbitrary values. Let $L$
be the set of $l$ for which $\sum_iQ^\infty_-(il)>0$. Then
$Q^\infty(ij)=\sum_{l\in L}Q^\infty_-(il)Q^\infty_+(lj)$ and the
$Q^t$ converge to $Q^\infty$. This proves the theorem.
\hfill$\square$
\begin{remark}\label{remark:strange}
{\em Theorem~\ref{thm:conv} says nothing of the convergence of the
$Q_+^t(lj)$ for those $l$ where $\sum_iQ^\infty_-(il)=0$. But
their behavior is uninteresting from a factorization point of
view. Indeed, since the $l$-th column of $Q^\infty_-$ is zero, the
values of the $l$-th row of $Q^\infty_+$ are not relevant, since
they don't appear in the product $Q^\infty_-Q^\infty_+$. As a
matter of fact, we now deal with an approximate nonnegative
factorization with a lower inner size. See also
remark~\ref{remark:lessk}. }
\end{remark}
\noindent In the next theorem we characterize the properties of
the fixed points of the algorithm. Recall from
section~\ref{sec:problem} that the objective function has no local
maxima in the interior of the domain.

\begin{thm}\label{thm:fixed}
If $(Q_-,Q_+)$ is a limit point of the
algorithm~(\ref{eq:q-bar2}),~(\ref{eq:q+bar2}) in the interior of
the domain, then it is a stationary point of the objective
function $D$. If $(Q_-,Q_+)$ is a limit point  on the boundary of
the domain corresponding to an approximate factorization where
none of the columns of $Q_-$ is zero ($\sum_iQ_-(il)>0$ for all
$l$), then all partial derivatives $\frac{\partial D}{\partial
Q_-(il)}$ and $\frac{\partial D}{\partial Q_+(lj)}$ are
nonnegative.
\end{thm}
{\bf Proof.} By computing the first order partial derivatives of
the objective function, using the middle term of
equation~(\ref{eq:divnew}), we can rewrite the update
equations~(\ref{eq:q-bar2}),~(\ref{eq:q+bar2}) as
\begin{equation}\label{eq:qn-}
Q_-^{t+1}(il) = Q_-^t (il) \left(-\frac{\partial D^t}{\partial
Q_-(il)} + 1 \right)
\end{equation}
and
\begin{equation}\label{eq:qn+}
Q_+^{t+1}(lj) \left(\sum_{i}Q_-^{t+1}(il) \right) = Q_+^t(lj)
\left(-\frac{\partial D^t}{\partial Q_+(lj)}+\sum_i
Q_-^t(il)\right).
\end{equation}
where $\frac{\partial D^t}{\partial Q_-(il)}$ stands for the
partial derivative $\frac{\partial D}{\partial Q_-(il)}$ evaluated
at $(Q_-^t,Q_+^t)$ and likewise for $\frac{\partial D^t}{\partial
Q_+(lj)}$.

Let $(Q_-,Q_+)$ be a limit point of the algorithm.
Equations~(\ref{eq:qn-}) and~(\ref{eq:qn+})  become
\[
Q_-(il)=Q_-{il} \left(-\frac{\partial D}{\partial Q_-(il)} +
1\right)
\]
\[
Q_+(lj)\left(\sum_{i}Q_-(il)\right)=Q_+(lj)\left(-\frac{\partial
D}{\partial Q_+(lj)}+\sum_i Q_-(il)\right).
\]
It follows that we then have the relations
\[
Q_-(il)\,\,\frac{\partial D}{\partial Q_-(il)}=0
\] and
\[ Q_+(lj)\,\,\frac{\partial D}{\partial Q_+(lj)}=0. \]
We first consider $Q_-$. Suppose that for some $i$ and $l$ we have
$Q_-(il)>0$, then necessarily $\frac{\partial D}{\partial
Q_-(il)}=0$. Suppose now that for some $i,l$ we have $Q_-(il)=0$
and that $\frac{\partial D}{\partial Q_-(il)}<0$. Of course, by
continuity, this partial derivative will be negative in a
sufficiently small neighborhood of this limit point. Since we deal
with a limit point of the algorithm, we must have infinitely often
for the iterates that $Q_-^{t+1}(il)<Q_-^t(il)$.
From~(\ref{eq:qn-}) we then conclude that in these points we have
$\frac{\partial D}{\partial Q_-(il)}>0$. Clearly, this contradicts
our assumption of a negative partial derivative, since eventually
the iterates will be in the small neighborhood of the limit point,
where the partial derivative is positive. Hence, we conclude that
$\frac{\partial D}{\partial Q_-(il)}\geq 0$,  if $Q_-(il)=0$. The
proof of the companion statement for the $Q_+(lj)$ is similar. If
$Q_+(lj)>0$, the corresponding partial derivative is zero. Let $l$
be such that
 $Q_+(lj)=0$ and suppose that
we have that $\frac{\partial D}{\partial Q_+(lj)}<0$. If we run
the algorithm, then $\frac{\partial D^t}{\partial
Q_+(lj)}/\sum_iQ^{t+1}_-(il)$ converges to a negative limit,
 whereas $\sum_iQ^{t}_-(il)/\sum_iQ^{t+1}_-(il)$ converges to one.
 Hence there is $\eta>0$ such that eventually
$\frac{\partial D^t}{\partial
Q_+(lj)}/\sum_iQ^{t+1}_-(il)<-2\eta/3$ and
$\sum_iQ^{t}_-(il)/\sum_iQ^{t+1}_-(il)>1- \eta/3$. Hence
eventually we would have, see~(\ref{eq:qn+}),
\[
Q_+^{t+1}(lj) -Q_+^{t}(lj) = Q_+^t(lj) \left(-\frac{\frac{\partial
D^t}{\partial Q_+(lj)}}{\sum_i Q_-^{t+1}(il)}+\frac{\sum_i
Q_-^t(il)}{\sum_i Q_-^{t+1}(il)} -1\right)>\eta/3,
\]
which contradicts convergence of $Q_+^{t}(lj)$ to zero.
\hfill$\square$

\begin{remark}
{\em If it happens that  a limit point $Q_-$ has a zero $l$-th
column, then it can easily be shown that the partial derivatives
$\frac{\partial D}{\partial Q_+(lj)}$ of $D$ are zero. Nothing can
be said of the values of the partial derivatives $\frac{\partial
D}{\partial Q_-(il)}$ for such $l$. But, see also
remark~\ref{remark:strange}, this case can be reduced to one with
a lower inner size factorization, for which the assertion of
theorem~\ref{thm:fixed} is valid. }
\end{remark}
\rm
\begin{cor}
The limit points of the algorithm with $\sum_iQ_-(il)>0$ for all
$l$ are all Kuhn-Tucker points for minimization of $D$ under the
inequality constraints $Q_-\geq 0$ and $Q_+\geq 0$.
\end{cor}
{\bf Proof.} Consider the Lagrange function $L$ defined by
\[
L(Q_-,Q_+)=D(P||Q_-Q_+)-\lambda\cdot Q_--\mu\cdot Q_+, \] where
for instance the inner product $\lambda\cdot Q_-$ is to be read as
$\sum_{il}\lambda_{il}Q_-(il)$ for $\lambda_{il}\in\mathbb{R}$.
Let us focus on a partial derivative $\frac{\partial L}{\partial
Q_-(il)}$ in a fixed point of the algorithm. The treatment of the
other partial derivatives is similar. From the proof of
theorem~\ref{thm:fixed} we know that in a fixed point we have
$Q_-(il)\frac{\partial D}{\partial Q_-(il)}=0$. Suppose that
$Q_-(il)>0$, then $\frac{\partial D}{\partial Q_-(il)}=0$ and the
Kuhn-Tucker conditions for this variable are satisfied with
$\lambda_{il}=0$. If $Q_-(il)=0$, then we know from
theorem~\ref{thm:fixed} that $\frac{\partial D}{\partial
Q_-(il)}\geq 0$. By taking $\lambda_{il}=\frac{\partial
D}{\partial Q_-(il)}\geq 0$, we see that also here the Kuhn-Tucker
conditions are satisfied. \hfill$\square$
\begin{remark}
{\em Wu~\cite{wu} has a number of theorems that characterize the
limit points of the closely related EM algorithm, or generalized
EM algorithm. These are all consequence of a general convergence
result in Zangwill~\cite{zangwill}. The difference of our results
with his is, that we also {\em have to} consider possible limit
points on the boundary, whereas Wu's results are based on the
assumption that all limit points lie in the interior of the
domain.}
\end{remark}\rm

\section{Relation with other minimization problems}\label{section:otherprob}

Other data analysis methods proposed in the literature enforce
some form of positivity constraint and it is useful to investigate
the connection between NMF and these methods. An interesting
example is the so called Archetypal Analysis (AA) technique
\cite{cb1994}. Assigned a matrix $X \in \R^{m\times n}$ and an
integer $k$, the AA problem is to find, in the convex hull of the
columns of $X$, a set of $k$ vectors whose convex combinations can
optimally represent $X$. To understand the relation between NMF
and AA we choose the $L_2$ criterion for both problems. For any
matrix $A$ and positive definite matrix $\Sigma$ define
$||A||_\Sigma = (\rm{tr} (A^T \Sigma A))^{1/2} $. Denote $||A||_I
= ||A||$. The solution of the NMF problem is then
$$
(W, H) = \arg \min_{W, H} || V - WH ||
$$
where the minimization is constrained to the proper set of
matrices. The solution to the AA problem is given by the pair of
column stochastic matrices $(A, B)$ of respective sizes $k \times
n$ and $m \times k$ such that $||X - XBA||$ is minimized (the
constraint to column stochastic matrices is imposed by the
convexity). Since $||X - XBA|| = ||I - BA||_{X^TX}$ the solution
of the AA problem is
$$
(A, B) = \arg \min_{A, B} || I - BA ||_{X^TX}.
$$
AA and NMF can therefore be viewed as special cases of a more
general problem which can be stated as follows. Given any matrix
$P \in \R_+^{m\times n}$, any positive definite matrix $\Sigma$,
and any integer $k$, find the best nonnegative factorization $P
\approx Q_1Q_2$ (with $Q_1 \in\R_+^{m\times k}, \,\, Q_2
\in\R_+^{k\times n}$) in the $L_2$ sense, {\it i.e.}
$$
(Q_1, Q_2) =  \arg \min_{Q_1, Q_2} ||P - Q_1Q_2||_\Sigma.
$$
{\bf Acknowledgement.} An anonymous referee is gratefully
acknowledged for helping us to improve the quality of the
presentation and for suggesting to us to investigate the boundary
behavior of the algorithm, similar to what has been reported
in~\cite{laa2004}. \rm

\end{document}